\documentclass[12pt,twoside,a4paper]{amsart}
\usepackage{amssymb}
\date{\today}

\title[On the Brian\c{c}on-Skoda theorem on a singular variety]
{On the Brian\c{c}on-Skoda theorem on a singular variety}

\author{Mats Andersson \& H{\aa}kan Samuelsson \& Jacob Sznajdman}
\thanks{The authors would like to thank the Institut Mittag-Leffler
(Djursholm, Sweden) where parts of this work were carried out.
The first author was
  partially supported by the Swedish 
  Research Council.}
\address{Department of Mathematics, Chalmers University of Technology and 
University of Gothenburg, S-412 96 G\"{O}TEBORG, SWEDEN}
\email{matsa@math.chalmers.se, hasam@math.chalmers.se, sznajdma@math.chalmers.se}

\def\ann{{\rm ann}}
\def\deg{\text{deg}\,}

\def\w{\wedge}

\def\dbar{\bar\partial}

\def\C{{\mathbb C}}
\def\w{{\wedge}}

\def\I{{\mathcal I}}

\def\Hom{{\rm Hom\, }}
\def\codim{{\rm codim\,}}

\def\Re{{\rm Re\,  }}

\def\ann{{\rm ann\,}}

\def\J{{\mathcal J}}
\def\I{{\mathcal I}}
\def\nbh{neighborhood }

\def\be{\begin{equation}}
\def\ee{\end{equation}}
\def\Ok{{\mathcal O}}

\newtheorem{thm}{Theorem}[section]

\newtheorem{cor}[thm]{Corollary}
\newtheorem{prop}[thm]{Proposition}

\theoremstyle{definition}

\theoremstyle{remark}

\newtheorem{preremark}{Remark}
\newtheorem{preex}{Example}

\newenvironment{remark}{\begin{preremark}}{\qed\end{preremark}}
\newenvironment{ex}{\begin{preex}}{\qed\end{preex}}

\numberwithin{equation}{section}

\begin{document}

\begin{abstract}
Let $Z$ be a  germ  of a reduced  analytic space of pure dimension.
We provide an analytic proof of the
uniform Brian\c con-Skoda theorem for the local ring $\Ok_Z$;
a result which was previously proved by Huneke by
algebraic methods.
For ideals with few generators we also get much
sharper  results.
\end{abstract}

\maketitle

\section{Introduction}\label{introduction}
Let $\mathfrak{a}=(a)=(a_1,\ldots, a_m)$ be an ideal in the local ring $\Ok=\Ok_0$ of 
holomorphic functions
at $0\in\C^d$ and let $|\mathfrak{a}|=\sum_j|a_j|$. Up to constants, this function is
independent of the choice of generators of $\mathfrak{a}$. 
In \cite{BS},  Brian\c{c}on and Skoda proved:

\smallskip
\noindent {\it If $\phi\in\Ok$ and 
\begin{equation}\label{bsvillkor}
|\phi|\le C|\mathfrak{a}|^{\min(m,d)+\ell-1},  \quad  \ell=1,2,3,\ldots,
\end{equation}
then $\phi\in \mathfrak{a}^\ell$.}

\smallskip
If $m\le d$, then the statement  follows directly from Skoda's $L^2$-estimate in
\cite{Sko};  if $m>d$ one  uses that there is an ideal $\mathfrak{b}\subset \mathfrak{a}$
such that $|\mathfrak{a}|\sim |\mathfrak{b}|$, a so-called reduction of $\mathfrak{a}$,  with
$n$ generators.
\smallskip

If  $\mathfrak{b}$ is any ideal in $\Ok$ then 
$|\phi|\le C|\mathfrak{b}|$ if (and in fact only if)  $\phi$ is in theintegral closure 
$\overline{\mathfrak{b}}$.
Therefore, the statement implies (is equivalent to) the inclusion
\begin{equation}\label{smaka}
\overline{\mathfrak{a}^{\min(m,d)+\ell-1}}\subset \mathfrak{a}^{\ell}.
\end{equation}
This  is a notable example of a purely algebraic theorem  that was first proved
by transcendental methods. It took several years before
algebraic proofs appeared,  \cite{LT} and 
\cite{LS}. 
In \cite{BGVY} there is a proof by integral formulas  and residue theory.

\smallskip

Assume now that  $Z$ is  a  germ of an analytic space of pure dimension $d$
and let $\Ok_Z$ be its  structure ring of germs of 
(strongly) holomorphic functions.
It is non-regular if (and only if) $Z$ is non-regular.
It is easy to see that the usual Brian\c con-Skoda theorem  cannot hold in general 
in the non-regular case,  not even for $m=1$, see Example~\ref{pucko} below.
However, 
Huneke proved in \cite{Hu}  that there is a number $\mu$ only
depending on $Z$ such that for any ideal $\mathfrak{a}\subset \mathcal{O}_Z$, and integer
$\ell\ge 1$,  
\begin{equation}\label{hun}
\overline{\mathfrak{a}^{\mu+\ell-1}}\subset \mathfrak{a}^{\ell}.
\end{equation}
Huneke's proof is completely algebraic (and holds for some more general rings as well),
so it is natural to look for an analytic proof.
In this paper we give a proof  by means of  residue calculus,
and the membership can be realized by an  integral formula 
on  $Z$. 
A problem of general interest, see, e.g., p.\ 657 in \cite{HyVi} and Remark 4.14 in \cite{Hu}, 
is to estimate the Brian\c{c}on-Skoda 
number, $\mu$, in Huneke's theorem in terms of invariants of the ring.
Our proof relates $\mu$ to the complexity of a free resolution of $\mathcal{O}_Z$.
We have also a sharper statement in case
$\mathfrak{a}$ has  ``few'' generators, and the zero set, $Z^{\mathfrak{a}}$, of the ideal does not
overlap the singular set of $Z$ ``too much''.  
To formulate this we first have to introduce  certain (germs of) subvarieties, 
$Z^r$, associated with $Z$:

To begin with we  choose an embedding of $Z$ and 
consider it as a subvariety at, say,
the origin of $\C^n$ for some $n$. If $\I$ is the corresponding radical ideal
in $\Ok=\Ok_{\C^n,0}$, then $\Ok_Z=\Ok/\I$.  
Let 
\begin{equation}\label{upplosning}
0\to \Ok(E_N)\stackrel{f_N}{\longrightarrow}\ldots
\stackrel{f_3}{\longrightarrow} \Ok(E_2)\stackrel{f_2}{\longrightarrow}
\Ok(E_1)\stackrel{f_1}{\longrightarrow}\Ok(E_0)
\end{equation}
be a free resolution of  $\Ok/\I$.  Here  $E_k$ are trivial vector bundles
and $E_0$ is a  trivial line bundle.
Thus $f_k$ are just   holomorphic matrices in a neighborhood of $0$.
We let 
$Z_k$ be the set of  points $x$ such that $f_k(x)$ does not have optimal rank.
These varieties  are,  see, \cite{Eis} Ch.~20,
independent of the choice of resolution, and  we have the inclusions
$$
\cdots\subset Z_{p+2}\subset Z_{p+1}\subset Z_{sing}\subset Z_p= \dots =Z_1=Z,
$$
where $p=n-d$.  
Now let     
\begin{equation}\label{hoppet}
Z^0=Z_{sing}, \quad  Z^r=Z_{p+r}, \ r>0.
\end{equation}
Since any  two minimal embeddings are equivalent, and any embedding factors
in a simple way over a minimal embedding, one can verify that these subsets
$Z^r$ are intrinsic subvarieties of the analytic space $Z$, that  reflect  the
degree of complexity of $Z$.   
To begin with, since $Z$ has pure dimension (Corollary 20.14 in \cite{Eis}),
$$
\codim Z^r\ge r+1,\quad r>0.
$$
Moreover, $Z^r=\emptyset$ for $r > d-\nu$ if and only if the 
depth of the ring $\mathcal{O}_Z$ is  at least $\nu$.
In particular,
$Z^r=\emptyset$ for $r>0$ if and only if $Z$ (i.e., $\Ok_Z$)  is
Cohen-Macaulay.

\begin{thm}\label{boll}
Let $Z$ be a germ of an  analytic space  of pure dimension.

\noindent (i)\  There is a natural number $\mu$, only depending on $Z$,
such that for any ideal $\mathfrak{a}=(a_1,\ldots,a_m)$ in $\Ok_Z$ and $\phi\in\mathcal{O}_Z$, 
\begin{equation}\label{bolla}
|\phi|\le C|\mathfrak{a}|^{\mu+\ell-1}
\end{equation}
implies that $\phi\in \mathfrak{a}^{\ell}$.

\smallskip
\noindent (ii) \   If  for a given ideal $\mathfrak{a}=(a_1,\ldots,a_m)$ 
\begin{equation}\label{bolla2}
\quad {\rm codim}\,(Z^r\cap Z^{\mathfrak{a}})\ge m+1+r,\quad   r\ge 0,
\end{equation}
then for any $\phi\in\Ok_Z$,
\begin{equation}\label{bolla3}
|\phi|\le C|\mathfrak{a}|^{m+\ell-1}
\end{equation}
implies that $\phi\in \mathfrak{a}^{\ell}$. 
\end{thm}

Huneke's theorem \eqref{hun}  follows immediately from part (i) of 
Theorem~\ref{boll}, since  even in the non-regular case
$\phi\in \overline{(b)}$ immediately implies that
$|\phi|\le C|b|$. The less obvious implication 
$|\phi|\le C|b| \Rightarrow \phi\in \overline{(b)}$ also holds, see, 
e.g., \cite{LeTe}, and so Theorem~\ref{boll} (i) is in fact equivalent 
to Huneke's theorem.

\begin{ex}\label{pucko}
If $Z$ is the zero set of $z^p-w^2$ in $\C^2$, where $p>2$ is a prime, 
then $|w|\le|z|^{[p/2]}$ on $Z$, but $w$ is not in $(z)$. However,
if $|\phi|\le C|z|^{(p+1)/2}$, then $\phi\in (z)$, i.e., $\phi/z$ is strongly
holomorphic on $Z$. 
\end{ex}

\begin{remark}
The important point in Huneke's  theorem  is the uniformity 
in $\mathfrak{a}$ and $\ell$. Notice that \eqref{hun} implies the slightly weaker statement 
\begin{equation}\label{hun2}
\overline{\mathfrak{a}}^{\mu +\ell -1}\subset \mathfrak{a}^{\ell}.
\end{equation}
It is quite easy to prove such an inclusion for fixed $\mathfrak{a}$ and $\ell$.
In fact, assume that $Z$ is a germ of a subvariety in $\C^n$
and choose a tuple $f$ such that  $Z=\{f=0\}$.  Let $A=(A_1,\ldots,A_m)$ 
and $\Phi$ denote
fixed representatives in $\Ok_{\C^n}$ of $\mathfrak{a}=(a_1,\ldots,a_m)$ and 
$\phi\in \overline{\mathfrak{a}}$. Then
\begin{equation*}
|\Phi(z)|\le C d(z,Z^{\mathfrak{a}}\cap Z)\le C'(|A|+|f|)^{1/M}
\end{equation*}
for some $M$ by  Lojasiewicz' inequality, and hence $\Phi^{Mn}$ is in
the ideal $(A)+(f)$ by the usual Brian\c{c}on-Skoda theorem
in the ambient space.  Thus $\phi^{Mn}\in \mathfrak{a}$ 
and 
therefore $\phi^{Mn\ell}\in\mathfrak{a}^\ell$. Thus 
$\overline{\mathfrak{a}}^{Mn\ell}\subset \mathfrak{a}^\ell$.
\end{remark}

From  Theorem~\ref{boll} (ii) we get:

\begin{cor}
If
\begin{equation}\label{bolla4}
{\rm codim}\,Z^r\ge m+1+r,\quad r\ge 0, 
\end{equation}
then  \eqref{bolla3} implies that $\phi\in\mathfrak{a}^\ell$ 
for any $\mathfrak{a}$ with $m$ generators.
\end{cor}

Assume that \eqref{bolla4} holds for $m=1$. The conclusion for $\ell =1$ then is that
each weakly holomorphic function is indeed holomorphic, i.e.,
$Z$ (or equivalently $\Ok_Z$) is normal. In fact,
if $\phi$ is weakly holomorphic, i.e., holomorphic on $Z_{reg}$ and
locally bounded, then it is meromorphic, so $\phi=g/h$ for some
$g,h\in\Ok_Z$. The boundedness means that $|g|\le C|h|$ and by the
corollary thus $\phi$ is in $\Ok_Z$. 
One can check that \eqref{bolla4} with $m=1$ is equivalent to
Serre's condition for normality of the local ring $\Ok_Z$ and therefore
both necessary and sufficient.


\smallskip

The basic tool in our proof is the residue calculus 
developed in \cite{A0}, 
\cite{AW1},  and \cite{AW2},  and we recall the necessary material
in  Section~\ref{utter}. Given an ideal sheaf $\J$ one 
can associate a current $R$ such that a holomorphic function
$\phi$ is in $\J$ as soon as $\phi R=0$. We use  such
a current $R^{\mathfrak{a},\ell}$ associated with the ideal $\mathfrak{a}^\ell$.
For $\ell=1$ it is the current of Bochner-Martinelli type
from \cite{A0}, whereas for $\ell>1$ we use a 
variant  from \cite{A2}. Since we are to prove the membership
on $Z$ rather than on some ambient space, thinking of $Z$ as embedded in some
$\C^n$, we will also use a  current $R^Z$
associated to the  radical ideal $I$ of the embedding.  For the analysis
of this current we rely  on  results from \cite{AS}, 
described in Section~\ref{zres}. 
It turns out that one can form the ``product'' $R^{\mathfrak{a},\ell}\w R^Z$ such that
$\phi R^{\mathfrak{a},\ell}\w R^Z$ only depends on the values of $\phi$ on $Z$;
moreover,  if the hypotheses in Theorem~\ref{boll} are fulfilled
then it vanishes (Proposition~\ref{main}), 
which  in turn implies that $\phi$ belongs to the ideal $\mathfrak{a}$ modulo $I$.
In the last section we  present   an integral formula that
provides an explicit representation of the membership.

\section{Currents obtained from locally free complexes}\label{utter}

Let 
\begin{equation}\label{vcomplex}
0\to E_N\stackrel{f_N}{\longrightarrow} E_{N-1}\stackrel{f_{N-1}}{\longrightarrow}\cdots
\stackrel{f_2}{\longrightarrow}E_1\stackrel{f_1}{\longrightarrow} E_0\to 0
\end{equation}
be a generically exact holomorphic complex of  Hermitian vector bundles
over a  complex manifold $X$, say a \nbh of the origin in $\C^n$. We  
assume that $E_0$ is a trivial line bundle so that 
$\Ok(E_0)=\Ok$.  There is an associated complex, like 
\eqref{upplosning},  of (locally) free sheaves
of $\Ok$-modules,  and  we let $\J=f_1\Ok(E_1)\subset\Ok$ be
the ideal sheaf generated by (the entries in) $f_1$. 
Let $Z$ be the analytic set where  \eqref{vcomplex} is not pointwise
exact. In $X\setminus Z$ we let $\sigma_k$ be the section of
$\Hom(E_{k-1},E_k)$ that vanishes on the orthogonal complement
of the pointwise image of $f_k$ and is  the minimal left inverse
of $f_k$ on the image of $f_k$. If $E=\oplus E_k$,
$f=\oplus f_k$, and $\sigma=\oplus \sigma_k$, then 
$
\sigma f+f\sigma=I,
$
where $I$ is the identity on $E$.
Since $E_0$ is trivial  we identify $\Hom(E_0,E)$ with $E$. 
Following \cite{AW1}, in $X\setminus Z$ we define the form-valued sections 
\begin{equation}\label{startmal}
u=\sum_{k=1}^N u_k, \quad  u_k=(\dbar\sigma_k)\cdots(\dbar\sigma_2)\sigma_1,
\end{equation}
of $E$. 
If $\nabla_f=f-\dbar$ we have that $\nabla_f u=1$.
It turns out that $u$ has a current  extension $U$ to
$X$ as a principal value current: If $F$ is a tuple
of holomorphic functions such that $F=0$ on $Z$, then
$|F|^{2\lambda} u$ has a current-valued  analytic continuation
to $\Re \lambda>-\epsilon$ and $U$ is the value at $\lambda=0$.
Alternatively one can take a smooth approximand of the characteristic
function $\chi$ for $[1,\infty)$, and  let $\chi_\delta=\chi(|F|^2/\delta^2)$.
Then $U$ is the weak limit of $\chi_\delta u$ when $\delta\to 0$
(see, e.g., the proofs of Theorems 16 and 21 in \cite{Sam}).
In this paper the latter definition will be more convenient.
Clearly $\nabla_f$ also applies to currents, and 
\begin{equation}\label{bulla}
\nabla_f U=1-R,
\end{equation}
where $R$ is a residue current with support on $Z$;  more precisely
$R=\lim_{\delta\to 0} R^\delta$,
where
$$
R^\delta=R^\delta_0+R^\delta_1+R^\delta_2+ \cdots =
(1-\chi_\delta)  +\dbar\chi_\delta\w u_1+\dbar\chi_\delta\w u_2+\cdots;
$$
notice that  $R^\delta_k$ is an $E_k$-valued $(0,k)$-current.

A basic observation is that the annihilator sheaf, $\ann R$, of 
$R$ is contained in the sheaf $\J$, i.e.,
\begin{equation}\label{pulla}
\ann R\subset \J.
\end{equation}
In fact, if  $\phi\in\Ok$ and $\phi R=0$, then by \eqref{bulla}, 
$
\nabla_f(\phi U)=\phi-\phi R=\phi. 
$
By  solving a sequence of $\dbar$-equations, which is always possible
locally at least, we get a holomorphic solution $\psi\in\Ok(E_1)$
to $f_1\psi=\phi$, which means that $\phi$ is in the ideal $\J$.
One can also prove \eqref{pulla} by an integral formula
that gives an explicit realization of the membership of $\phi$ in $\J$, 
see Section~\ref{intrep}. 

In general the converse inclusion is not true. However,
if the associated sheaf complex
is exact, i.e., a resolution of $\Ok/\J$, then
indeed  $\ann R=\J$ (Theorem~1.1 in \cite{AW1}).

\begin{ex}\label{bsex}  
Let $a_1,\ldots,a_m$ be holomorphic functions in $X$.
Choose a nonsense basis $\{e_1,\ldots,e_m\}$ and  consider
$E_1=sp\{e_j\}$ as a trivial vector bundle of rank $m$, let $e_j^*$
be the dual basis, and  consider 
$a=a_1e_1^*+\cdots +a_me_m^*$ as a section of the dual bundle $E_1^*$.
If $E_k=\Lambda^k E_1$
we then get a complex \eqref{vcomplex},
the Koszul complex, 
with the mappings  $f_k$ as interior multiplication $\delta_a$ with $a$.
Following the recipe above (with the trivial metric on the $E_k$) we get, cf.,  \cite{AW1} Example~1,
the corresponding form 
\begin{eqnarray}\label{equ}
u^a=\sum_{k=1}^{m}
\frac{(\sum_{j=1}^m \bar{a}_je_j)\wedge (\sum_{j=1}^m \dbar \bar{a}_j\wedge e_j)^{k-1}}{|a|^{2k}}
\end{eqnarray}
outside $\{a=0\}$
and the associated residue current $R^a=\lim_{\delta\to 0} R^{a,\delta}$
where 
$R^{a,\delta}= (1-\chi_\delta)+
\dbar\chi_\delta\wedge u^a$ and $\chi_\delta=\chi(|a|^2/\delta^2)
$. 
This current of so-called Bochner-Martinelli type
was introduced already in \cite{PTY}, and its  relation to the
Koszul complex and division problems was noticed in  \cite{A0}.
Now \eqref{pulla} means that  
\begin{equation}\label{ain}
\ann R^a\subset (a).
\end{equation}
Except for the case when $a$ is a complete intersection,
in which case the Koszul complex provides a resolution of
$\Ok/(a)$, the   inclusion \eqref{ain} is strict, see \cite{W} and \cite{JW}.  
Nevertheless, the singularities of  $R^a$  
reflect the
characteristic varieties associated to the ideal, see \cite{JW} and \cite{AG},
which are  closely related to the integral closure
of powers of $(a)$,  and therefore
$R^a$   is well suited for the Brian\c{c}on-Skoda theorem. 

A slight modification of the Koszul complex, derived  from the so-called
Eagon-Northcott complex, with associated ideal sheaf $\J=(a)^\ell$, 
was introduced in \cite{A2}.
The associated form $u^{a,\ell}$ is a sum of terms like
$$
\frac{\bar a_{I_1}\cdots\bar a_{I_\ell}\dbar\bar a_{I_{\ell+1}}\w\ldots\w \dbar\bar a_{I_{k+\ell-1}}}
{|a|^{2(k+\ell-1)}},\quad k\le m,
$$
see the proof of Theorem~1.1 in \cite{A2} for 
a  precise description of $u^{a,\ell}$ and the corresponding residue current $R^{a,\ell}$.
It turns out that  $\phi$ annihilates   $R^{a,\ell}$
if \eqref{bsvillkor} holds,   and thus $\phi\in(a)^\ell$, so the
classical Brian\c con-Skoda theorem follows.
The most expedient way to  prove 
this annihilation  is to use a  resolution
of singularities where $a$ is principal. 
However, it is not really necessary to define 
the current $R^{a,\ell}$ in itself; 
it is actually enough to make sure that $\phi R^{a,\ell,\delta}\to 0$
when $\delta\to 0$, and this can be proved 
essentially  by integration by part in an ingenious way, thus providing
a proof of  the Brian\c{c}on-Skoda theorem by 
completely elementary means, see \cite{Sz}. 
\end{ex}

In  \cite{AW2} was  introduced the sheaf of  {\it pseudomeromorphic}
currents $\mathcal{PM}$.  For the definition,  
see \cite{AW2}. 
It is closed under $\dbar$ and  multiplication with
smooth forms.  In particular, the  currents $U$ and $R$ are pseudomeromorphic.
The following fact (Corollary~2.4 in \cite{AW2})  will be used
repeatedly.

\begin{prop}\label{dimer}
If $T\in\mathcal{PM}$ has bidegree $(r,k)$ and the support of $T$ is contained
in a variety of codimension strictly larger than $k$, then $T=0$.
\end{prop}

In particular, this means that if $Z$ (the variety where \eqref{vcomplex} is not pointwise
exact) has codimension $p$ then $R=R_p+R_{p+1}+\cdots$.

\smallskip

As mentioned in the introduction, we need to form products of currents associated to complexes.
Assume therefore that $(\Ok(E^g_{\bullet}), g_{\bullet})$ and $(\Ok(E^h_{\bullet}), h_{\bullet})$ 
are two complexes as above
and $\I$ and $\J$ are the corresponding ideal sheaves.
We can  define a   complex  \eqref{vcomplex}  with 
\begin{equation}\label{sven}
E_k=\bigoplus_{i+j=k} E^g_i\otimes E^h_j,
\end{equation}
and $f=g+h$, or more formally,
$f=g\otimes I_{E^h}+ I_{E^g}\otimes h$,
such that
\begin{equation}\label{hoppsan}
f(\xi\otimes\eta)=g\xi\otimes\eta +(-1)^{\deg\xi}\xi\otimes h\eta.
\end{equation}
Notice that $E_0=E_0^g\otimes E_0^h=\C$ and that
$f_1\Ok(E_1)=\I+\J$.
One can extend \eqref{hoppsan} to form-valued  or current-valued sections $\xi$ and $\eta$ 
and $\deg\xi$ then means  total degree. It is  natural to write
$\xi\w \eta$ rather than $\xi\otimes\eta$, and 
we  define $\eta\w\xi$  as $(-1)^{\deg \xi \deg\eta}\xi\w\eta$.
Notice that 
\begin{equation}\label{rregel}
\nabla_f (\xi\otimes\eta)=\nabla_g\xi\otimes\eta
+(-1)^{\deg \xi}\xi\otimes \nabla_h\eta.
\end{equation}
Let  $u^g$ and $u^h$ be the corresponding
$E^g$-valued and $E^h$-valued 
forms, cf.\ \eqref{startmal}.  
Then $u=u^h\w u^g$ is an  $E$-valued form
outside $Z^g\cup Z^h$. 
Following the proof of  Proposition~2.1 in \cite{AW2} we can define
$E$-valued pseudomeromorphic currents
$$
R^h\w  R^g=\lim_{\delta\to 0} R^{h,\delta}\w  R^g,  
\quad 
U^h\w R^g=\lim_{\delta\to 0} U^{h,\delta}\w R^g,
$$
where
$U^{h,\delta}=\chi_\delta u^{h}$ and $R^{h,\delta}=1-\chi_{\delta}+\dbar\chi_\delta \w  u^{h}$,
and $\chi_\delta=\chi(|H|^2/\delta^2)$ as before. 
The ``product''  $R^h\w R^g$ so defined is {\it not} equal to
$R^g\w R^h$ in general. It is also understood here that
$H$ only vanishes where it has to, i.e., on the set where
the complex $(E^h_{\bullet}, h_{\bullet})$ is not pointwise exact. If we use 
an $H$ that vanishes on a larger set, the result will be
affected. It is worth to point out that  a certain component 
$R^h_k\w R^g$ may be nonzero even if $R^h_k$ itself vanishes.

\begin{prop}\label{bellman}
With the notation above we have that
\begin{equation}\label{armar}
\nabla_f (U^g+ U^h\w R^g)=1-R^h\w R^g.
\end{equation}
Moreover,  $\phi R^h\w R^g=0$ implies that
$\phi\in\I+\J$.
\end{prop}

\begin{proof}
Recall that $\nabla_h U^{h,\delta}=1-R^{h,\delta}$, $\nabla_g U^g=1-R^g$ and 
$\nabla_gR^g=0$. 
Therefore,
$$
\nabla_f (U^g+ U^{h,\delta}\w R^g)=1-R^g+(1-R^{h,\delta})\w R^g=
1-R^{h,\delta}\w R^g.
$$
Taking limits, we get \eqref{armar}. The second statement now follows
in the same way as \eqref{pulla} above.
\end{proof}



\section{The residue  current associated to the variety $Z$}\label{zres}

Consider a subvariety $Z$ of a \nbh of the origin in $\C^n$ with
radical ideal sheaf $\I$ and let \eqref{upplosning}
be a resolution of $\Ok/\I$.  Let $R^Z$ be the associated
residue current obtained as in the previous section.
We then know that $R^Z$ has support on $Z$ and that
$\ann R^Z=\I$.
Outside the set $Z_k$, cf., Section~\ref{introduction},  the mapping $f_k$ has constant
rank, and hence $\sigma_k$ is smooth there. Outside
$Z_k$ we therefore have that 
\begin{equation}\label{strumpa}
R^Z_{k+1}=\alpha_{k+1} R_k^Z
\end{equation}
where $\alpha_{k+1}=\dbar\sigma_{k+1}$ 
is a smooth $\mbox{Hom}(E_k,E_{k+1})$-valued $(0,1)$-form,
cf., \eqref{startmal}.


Locally on $Z_{reg}$, the current 
$R^Z$ is essentially the integration current $[Z]$.
We have the following more precise statement that
gives a Dolbeault-Lelong-type representation, in the sense of \cite{jeb1},
of the current $R^Z$. 
Let  $\chi$ be  a  smooth regularization of the characteristic function of
$[1,\infty)$ and $p=\textrm{codim}\, Z$ as before.

\begin{prop}\label{DL}
For each given $x\in Z_{reg}$, there  is a hypersurface $\{h=0\}$ in $Z$, avoiding $x$ but 
containing $Z_{sing}$ and intersecting
$Z$ properly,  and  $E_k$-valued $(n-p,k-p)$-forms 
$\beta_k$,  smooth outside $\{h=0\}$, such that
\begin{equation*}
R^Z_k .\, (dz\wedge \xi)=\lim_{\epsilon\to 0}
\int_{Z}\chi(|h|/\epsilon)\beta_k \wedge \xi,\,\,\,\,\,\, \xi \in \mathcal{D}_{0,n-k}(X),
\end{equation*}
for $p\leq k \leq n$. 
Moreover, in a suitable resolution $\pi\colon \tilde Z\to Z$
the forms $\beta_k$ locally have  the form $\alpha_k/m_k$, where $\alpha_k$ are  smooth
and $m_k$ are   monomials.
\end{prop}

Here, $dz =dz_1\wedge \cdots \wedge dz_n$.

\begin{proof}
Following  Section~5 in \cite{AS} (the proof of Proposition~2.2)
one can find, for each given $x\in Z_{reg}$, a holomorphic function $h$ such that $h(x)\neq 0$ 
and $h$ does not vanish identically on any component of $Z_{reg}$. Moreover, for $k\ge p$,
$$
R_k^Z=\gamma_k\lrcorner [Z],
$$
where $\gamma_k$ is an  $E_k$-valued and $(0,k-p)$-form-valued
$(p,0)$-vector field that is smooth outside $\{h=0\}$. 
Let $\xi$ be a test form of bidegree $(0,n-k)$.
The current $R^Z$  has the so-called standard extension property, SEP, see \cite{AW2}
Section~5, which means that 
$$
R^Z_k. (\xi\w dz)=\lim_{\epsilon\to 0}\int\chi(|h|/\epsilon)\gamma_k\lrcorner[Z]\w \xi\w dz=
\pm 
\lim_{\epsilon\to 0}\int_Z\chi(|h|/\epsilon)\xi\w \gamma_k\lrcorner dz.
$$
Thus we can take $\beta_k=\pm  \gamma_k\lrcorner dz$.

More precisely,  according to the  last paragraph of Section~5 in \cite{AS}, 
$\gamma_p$ is a 
meromorphic  $(p,0)$-field (with poles where $h=0$)
composed by the  orthogonal projection
of $E_p$ onto the orthogonal complement in $E_p$ of the pointwise image of $f_{p+1}$.
This projection is given by
$$
I_{E_p}-f_{p+1}\sigma_{p+1}.
$$
Furthermore, cf., \eqref{strumpa},
$$
\gamma_k=(\dbar\sigma_k)\cdots(\dbar\sigma_{p+1})\gamma_p
$$
for $k>p$. 
Now choose a resolution of singularities
$\tilde Z\to Z$ such that for each $k$ the
the determinant ideal of $f_k$ is principal. On $\tilde Z$, then
each  $\sigma_k$  (locally) is 
a smooth form over a monimial, see Section~2 in \cite{AW1}, and thus
$\beta_k=\gamma_k\lrcorner dz$  has this form as well. 
\end{proof}

We can choose the resolution of singularities $\tilde{Z} \to Z$ 
so that also $\tilde h=\pi^*h$ is a monomial.
By a partition of unity  it follows that 
$R^Z_k.(dz\w\xi)$ is a finite sum of terms like
\begin{equation}\label{armada}
\lim_{\epsilon\to 0} \int_s\chi(|\tilde h|/\epsilon)
\frac{ds_1\w\ldots\w ds_{\nu}}{s_1^{\alpha_1+1}\cdots s_{\nu}^{\alpha_{\nu}+1}}
\w\tilde\xi\w\psi,
\end{equation}
where $s_1,\ldots, s_{n-p}$ are local holomorphic coordinates and $\nu\le n-p$, 
$\tilde\xi=\pi^*\xi$, and $\psi$ is a smooth form with compact support.
It is easily checked that this limit  is the
tensor product of the one-variable principal value currents
$ds_i/s_i^{\alpha_i+1}$,
$1\le j\le\nu$, 
acting on $\tilde\xi\w\psi$. Therefore \eqref{armada}  is equal
to (a constant times)
\begin{equation}\label{armada2}
\int
\frac{ds_1\w\ldots\w ds_{\nu}}{s_1\cdots s_{\nu}}\w\partial^\alpha_s(\tilde \xi\w\psi),
\end{equation}
if $\partial^\alpha_s=\partial^{\alpha_1}_{s_1}\cdots \partial^{\alpha_{\nu}}_{s_\nu}$.

\section{Proof of Theorem \ref{boll}}\label{pf}

To prove Theorem~\ref{boll} we are going to apply the idea
in Example~\ref{bsex} but performed on $Z$. 
To this end we assume that $Z$ is embedded in $\C^n$ and we let
$R^Z$ be the current introduced in the previous section.
Let $\mathfrak{a}=(a)$ be the ideal in $\Ok_Z$ and suppose for the moment
that $a$ also denotes representatives in $\Ok$ of the generators.
If $R^{a,\ell}=\lim_{\delta\to 0} R^{a,\ell,\delta}$ denotes
the current from  Example~\ref{bsex} we can form,
cf., the end of Section~\ref{utter}, the product
$$
R^{a,\ell}\w R^Z=\lim_{\delta\to 0}R^{a,\ell,\delta}\w R^Z.
$$
Since $R^Z$ annihilates $\I$ it follows that $R^{a,\ell}\w R^Z$ only
depends on $\mathfrak{a}\subset\Ok_Z$.  For the same reason,
$\phi R^{a,\ell}\w R^Z$
is well-defined for $\phi\in\Ok_Z$.  We know from Proposition~\ref{bellman}
that
$\phi$ belongs to $\mathfrak{a}$ if it annihilates this current, and thus
Theorem~\ref{boll} follows from the following proposition.

\begin{prop}\label{main}
If the hypotheses of Theorem~\ref{boll} are fulfilled 
i.e., either \eqref{bolla}, or \eqref{bolla3} together with the geometric
conditions \eqref{bolla2}, then
$\phi R^{a,\ell}\w R^Z=0$.  
\end{prop}

\begin{remark}
It is natural to try to use the Lelong current $[Z]$ rather than $R^Z$.
There is, see  \cite{AB} Example~1,  a holomorphic $E_p$-valued
form $\xi$ such that $[Z]=\xi\cdot R^Z_p$. Thus the hypotheses
in Theorem~\ref{boll} imply that
$\phi R^a\w[Z]=0$. However, this in turn does not imply that
$\phi$ is in $(a)$. 
In fact, if $m=1$ so that $a$ is just one function, then 
$$
0=\phi R^a\w [Z]=\phi\dbar\frac{1}{a}\w[Z], 
$$ 
and this means that $\phi/a$ is in  $\omega_Z^0$ introduced by
Barlet, see, e.g., \cite{HP},
and this class  is wider than $\Ok_Z$ in general. 
\end{remark}

\begin{proof}[Proof of Proposition~\ref{main}]

We first  assume that \eqref{bolla2} and \eqref{bolla3} hold.
Considering $\phi R^{a,\ell}$ as an intrinsic  current on
the submanifold $Z_{reg}$ (cf.\ the beginning of this section)
it follows from the residue proof of the Brian\c{c}on-Skoda theorem
in the regular case that $\phi R^{a,\ell}$ must vanish on $Z_{reg}$
since \eqref{bolla3} holds. Thus, $\phi R^{a,\ell}\wedge [Z]$ vanishes on $Z_{reg}$
and so, 
in view of  Proposition~\ref{DL}, it follows that
the support of $\phi R^{a,\ell}\w R^Z$ is contained in $Z_{sing}$.
On the other hand it is readily verified that
$R^{a,\ell}\w R^Z$ must vanish if $a$ is nonvanishing. 
Thus the support of $\phi R^{a,\ell}\w R^Z$ is contained in $Z_{sing}\cap Z^a$.

The current  $R^{a,\ell}$ has (maximal)  bidegree $(0,m)$ and hence
$R^{a,\ell}\w R^Z_p$ has (maximal) bidegree $(0,m+p)$. Since it has support
on  $Z_{sing}\cap Z^a$ that has codimension $\geq p+m+1$ by \eqref{bolla2}, 
it follows that $\phi R^{a,\ell}\w R^Z_p=0$. 
Outside $Z_{p+1}$ we have that $R^Z_{p+1}=\alpha_{p+1} R^Z_p$ for a smooth form $\alpha_{p+1}$,
and hence
$$
\phi R^{a,\ell}\w R^Z_{p+1}=\phi R^{a,\ell}\w \alpha_{p+1} R^Z_p=\alpha_{p+1}\phi R^{a,\ell}\w R^Z_p =0
$$
there. Thus $\phi R^{a,\ell}\w R^Z_{p+1}$ has support on $Z_{p+1}\cap Z^a$, and again for degree
reasons we find  that $\phi R^{a,\ell}\w R^Z_{p+1}=0$. Continuing in this way
we can conclude that $\phi R^{a,\ell}\w R^Z=0$.

\smallskip
We now assume that \eqref{bolla} holds.
We have to prove that $R^Z.(dz\w\xi)\to 0$ when $\delta\to 0$, for
\begin{equation}\label{pelle}
\xi=\phi R^{a,\ell,\delta}\w\eta,
\end{equation}
with  test forms $\eta$ of bidegree $(0,*)$. 
In view of the comments after the proof of
Proposition~\ref{DL} it is enough to prove that 
each term \eqref{armada2} tends to zero if 
\eqref{bolla} holds and $\mu$ is large enough
(independently of $(a)$ and $\ell$). 
For this particular term we will see that we need $\mu\ge\mu_0$, where
\begin{equation}\label{val}
\mu_0=|\alpha|+2\min(m,n-p).
\end{equation}
For simplicity we omit all snakes from now on
and write $\phi$ rather than $\tilde \phi$ etc.
Moreover, we assume that $\ell=1$, the general case follows completely analogously.
Since $\tilde Z$ is smooth, by the usual Brian\c{c}on-Skoda theorem we have
that 
\begin{equation}\label{pelle1}
 \phi\in( a)^{|\alpha|+\min(m,n-p)+1}.
\end{equation}
Notice that
$$
R^{a,\delta}_{k}=\chi'(|a|^2/\delta^2)\w\frac{\dbar|a|^2}{\delta^2}\w u^a_k, \quad k>0,
$$
and thus $R^{a,\delta}_k$ is a sum of terms like
$$
\chi' \frac{\dbar\bar a_{I_1}\w\ldots\w\dbar\bar a_{I_k}}{\delta^2|a|^{2k}}\bar a a\w \omega
$$
for $|I|=k$, where in what follows
$a^r$ denotes a product of $r$ factors $a_i$, and similarly with $\bar a^r$,
and $\omega$ denotes a smooth form.
For degree reasons $k\le \nu=\min(m, n-p)$.
In view of \eqref{pelle1} therefore
$\phi R^{a,\delta}_{k}$ is a sum of terms like
$$
\chi' \frac{\dbar\bar a_{I_1}\w\ldots\w\dbar\bar a_{I_\nu}}{\delta^2|a|^{2\nu}}
\bar a a^{2+\nu+|\alpha|}\w \omega
$$
plus lower order terms.
A straight forward computation yields that
$\partial_s^\alpha (\phi R^{a,\delta}_{k})$ is a finite sum of terms like
$$
\chi^{(r+1)} \frac{\dbar\bar a_{I_1}\w\ldots\w\dbar\bar a_{I_\nu}}
{\delta^{2(r+1)} |a|^{2(\nu+|\gamma|-r)}} \bar a^{1+|\gamma|} a^{2+\nu+|\gamma|}\w \omega,
$$
where $\gamma\le\alpha$ and $r\le|\gamma|$, 
plus lower order terms.

We thus  have to see that each 
\begin{equation}\label{stolle}
\int_s 
\frac{ds_1\w\ldots\w ds_{\nu}}{s_1\cdots s_{\nu}}
\chi^{(r+1)} \frac{\dbar\bar a_{I_1}\w\ldots\w\dbar\bar a_{I_\nu}}
{\delta^{2(r+1)} |a|^{2(\nu+|\gamma|-r)}} \bar a^{1+|\gamma|} a^{2+\nu+|\gamma|}\w \omega
\end{equation}
tends to $0$ when $\delta\to 0$. 
After a suitable further resolution we may assume that locally
$a=a_0a'$ where $a_0$ is holomorphic and $a'$ is a non-vanishing
tuple.
Then 
$$
\dbar\bar a_{I_1}\w\ldots\w\dbar\bar a_{I_\nu}=\bar a_0^{\nu-1}\w \omega.
$$
Also notice that  the expression
\begin{equation}\label{stall}
\frac{ds_1\w\ldots\w ds_{\nu}}{s_1\cdots s_{\nu}}
\end{equation}
becomes a sum of similar expressions in this new  resolution.
Altogether we end up with a finite sum of terms like
$$
\int_s 
\frac{ds_1\w\ldots\w ds_{\nu}}{s_1\cdots s_{\nu}}
\chi^{(r+1)} (|a|^2/\delta^2)\w \Ok(1),
$$
and each such integral tends to zero by dominated convergence.

The term corresponding to 
$R^{a,\delta}_0=1-\chi(|a|^2/\delta^2)$ is handled in a similar but easier way.
\end{proof}

\section{Integral representation of the membership}\label{intrep}

Finally  we describe  how one can obtain an explicit integral
representation of the membership provided that the residue is annihilated.
The starting point is the formalism in \cite{A1} to generate integral
representations for holomorphic functions. Let $\delta_\eta$ denote
interior multiplication with the vector field
$$
2\pi i\sum_1^n (\zeta_j-z_j)\frac{\partial}{\partial\zeta_j}
$$
and let $\nabla_\eta=\delta_\eta-\dbar$.
A smooth form $g=g_0+g_1+\cdots +g_n$, where $g_k$ has bidegree $(k,k)$, is
called a {\it weight}  (with respect to $z$) if $\nabla_\eta g=0$ and $g_0(z,z)=1$.
Notice that the product of two weights is again a weight.

\begin{ex}\label{sopp}
Let $\chi$ be a cutoff function that is identically $1$ in a \nbh of the
closed unit ball, and let
$$
s=\frac{1}{2\pi i}\frac{\partial|\zeta|^2}{|\zeta|^2-\bar\zeta\cdot z}.
$$
Then $\nabla_\eta s=1-\dbar s$ and therefore 
$$
g=\chi -\dbar\chi\w [s+s\w\dbar s+\cdots +s\w(\dbar s)^{n-1}]
$$
is a weight with respect to $z$ for each $z$ in the ball, with compact support,
and  it depends holomorphically on $z$.
\end{ex}

If $g$ is a weight with compact support and $z$ is holomorphic on the support,
then
$$
\phi(z)=\int g\phi =\int g_n\phi.
$$

Now consider a complex like \eqref{vcomplex}  in Section~\ref{utter},
 defined in a \nbh
of the closed ball, and let $U^\delta$ and $R^\delta$ be the associated 
$E$-valued forms. 
One can find, see \cite{A11} Proposition~5.3, holomorphic $E^*_k$-valued
$(k,0)$-forms $H^0_k$ and $\textrm{Hom}(E_k,E_1)$-valued
$(k-1,0)$-forms $H^1_k$ such that $\delta_{\eta} H^0_k=H^0_{k-1} f_k(\zeta)-f_1(z) H^1_k$
and $H^j_j=\textrm{Id}_{E_j}$.
Using that $\nabla_{f}U^{\delta}=1-R^{\delta}$ one verifies that
$$
f_1(z)HU^\delta+HR^\delta = 1-\nabla_{\eta}(\sum H^0_kU_k^{\delta}),   
$$
where
$$ 
H U^{\delta}=\sum H^1_k U_k^\delta, \quad H R^\delta=\sum H^0_k R_k^\delta.
$$
It follows that $g^{\delta}:=f_1(z)HU^\delta+HR^\delta$ is a weight with respect to $z$.
If $g$ is, e.g., the weight from Example~\ref{sopp}  we thus get the 
representation
$$
\phi(z)=\int g^{\delta}\wedge g \phi=f_1(z)\int HU^\delta\w g \phi+\int HR^\delta\w g\phi. 
$$
Taking limits we obtain  the interpolation-division formula
\begin{equation}\label{pup}
\phi(z)=f_1(z)\int HU\w g \phi+\int HR\w g\phi.
\end{equation}
To be precise, the integrals here are the action of currents on smooth forms.
In particular, \eqref{pup} implies  that $\phi$ belongs to the ideal generated by $f_1$
if $\phi R=0$.

If we now choose as our complex the resolution of the sheaf $I=I_Z$, we get
the formula
$$
\phi(z)=\int g\w H^ZR^Z \phi, \quad z\in Z,
$$
for $\phi\in\Ok_Z$.
We then   replace $g$ by the weight
$g^{a,\ell,\delta}\w g$,
where
\begin{equation*}
g^{a,\ell,\delta}=a(z)^\ell\cdot  H^{a,\ell} U^{a,\ell,\delta}+ H^{a,\ell} R^{a,\ell,\delta};
\end{equation*} 
here  $a(z)^\ell$ denotes  the first mapping in the complex
associated with $(a)^\ell$, cf., Example~\ref{bsex},  so that its entries
are elements in the ideal $(a)^\ell$. 
We get 
\begin{eqnarray*}\label{eq2}
\phi(z) &=& a(z)^\ell\cdot\int_{\zeta} H^a U^{a,\ell,\delta}\w
H^ZR^Z \phi \wedge g  \\
& &
+\int_{\zeta} H^a R^{a,\ell,\delta}\w
H^ZR^Z\wedge g \phi. \nonumber
\end{eqnarray*}
If the hypotheses in Theorem~\ref{boll}
are fulfilled, since $H^Z$, $H^a$ and $g$ are smooth, the second integral tends
to zero when $\delta\to 0$, and the  first integral on the right hand side  
converges to an $E_1^{a,\ell}$-valued 
holomorphic function.
Thus we get the  explicit representation 
$$
\phi(z)=
 a(z)^\ell\cdot\int_{\zeta} H^a U^{a,\ell}\w
H^ZR^Z \phi \wedge g 
$$
of the membership.

\def\listing#1#2#3{{\sc #1}:\ {\it #2},\ #3.}


\begin{thebibliography}{9999}






\bibitem{A0}\listing{M.\ Andersson}
{Residue currents and ideals of holomorphic functions}
{Bull.\  Sci.\ Math., {\bf 128} (2004), 481--512}



\bibitem{A1}\listing{M.\ Andersson}
{Integral representation  with weights I}
{Math.\ Ann., {\bf 326} (2003), 1--18}




\bibitem{A11}\listing{M.\ Andersson}
{Integral representation with weights. II. Division and interpolation}
{Math.\ Z.,  {\bf 254}  (2006), 315--332}





\bibitem{A2} \listing{M.\ Andersson}
{Explicit versions of the Briancon-Skoda theorem with variations}
{Michigan Math.\ J., {\bf 54(2)} (2006), 361--373}




\bibitem{AB}\listing{M.\ Andersson}
{Coleff-Herrera currents, duality, and Noetherian operators}
{Preprint Gothenburg (2009),  \  available at arXiv:0902.3064} 







\bibitem{AS} \listing{M.\ Andersson, H.\ Samuelsson}
{Koppelman formulas and the  $\bar\partial$-equation on 
an analytic space}
{Preprint Mittag-Leffler  (2008), available at 	arXiv:0801.0710}











\bibitem{AW1} \listing{M.\ Andersson, E.\ Wulcan}
{Residue currents with prescribed annihilator ideals}
{Ann.\ Sci.\ \'{E}cole Norm.\ Sup.,  {\bf 40} (2007), 985--1007}



\bibitem{AW2} \listing{M.\ Andersson, E.\ Wulcan}
{Decomposition of residue currents}
{J.\ reine  angew.\ Math., (to appear), \  
available at 	arXiv:0710.2016 }







\bibitem{Astark} \listing{M.\ Andersson}
{A residue criterion for strong holomorphicity}
{Preprint Gothenburg 2007,
\  available at arXiv:0711.2863}









\bibitem{AG}\listing{M.\ Andersson, E.\ G\"otmark}
{Explicit representation of membership of polynomial ideals}
{Preprint Mittag-Leffler (2008),\  
 available at  arXiv:0806.2592 }






\bibitem{BGVY}\listing
{C.\ Berenstein, R.\ Gay, A.\ Vidras, 
A.\ Yger}
{Residue Currents and B\'ezout Identities}
{Birkh\"auser (1993)}




\bibitem{jeb1} \listing{J.-E. Bj\"{o}rk}
{Residue currents and $\mathcal{D}$-modules}
{The legacy of Niels Henrik Abel, 605--651, Springer, Berlin 2004}




\bibitem{BS}\listing{J.\  Brian\c{c}on, H.\ Skoda}
{Sur la cl\^oture int\'egrale d'un id\'eal de germes de fonctions 
holomorphes en un point de $\C\sp{n}$}
{ C.\ R.\ Acad.\ Sci.\ Paris S\'er. A 278 (1974),
949--951}





\bibitem{CH}\listing{N.r.\ Coleff, M.e.\ Herrera}
{Les courants r\'esiduels associ\'es \`a une forme m\'eromorphe}
{Lect. Notes in Math. {\bf 633},  Berlin-Heidelberg-New York (1978)}



\bibitem{Eis}\listing{D.\ Eisenbud}
{Commutative algebra.
With a view toward algebraic geometry}
{Graduate Texts in Mathematics, 150.
Springer-Verlag, New York, 1995}





\bibitem{HP}\listing{G.\ Henkin, M.\ Passare}
{Abelian differentials on singular varieties and variatiosn on a theorem
of Lie-Griffiths}
{Invent.\ Math., {\bf 135} (1999), 297--328}



\bibitem{Hu}\listing{C.\ Huneke}
{Uniform bounds in Noetherian rings}
{Invent.\ Math., {\bf 107}  (1992),  203--223}


\bibitem{HyVi}\listing{E.\ Hyry, O.\ Villamayor}
{A Brian\c{c}on-Skoda Theorem for Isolated Singularities}
{J.\ Algebra, {\bf 204}  (1998),  656--665}


\bibitem{JW}\listing{M.\ Jonsson, E.\ Wulcan}
{On Bochner-Martinelli residue currents and their annihilator ideals}
{arXiv:0811.0636 }{}{}{}{} 



\bibitem{LeTe}\listing{M.\ Lejeune-Jalabert, B.\ Teissier, J.-J.\ Risler}
{Cl\^{o}ture int\'{e}grale des id\'{e}aux et \'{e}quisingularit\'{e}}
{arXiv:0803.2369}




\bibitem{LS}
\listing{J.\ Lipman, A.\ Sathaye}
{Jacobian ideals and a theorem of Brian\c{c}on-Skoda}  
{Michigan Math.\ J., {\bf 28}  (1981), no. 2, 199--222}





\bibitem{LT}\listing{J.\ Lipman, B.\ Teissier}
{Pseudorational local rings and a theorem of
Brian\c{c}on-Skoda about integral closures of ideals}
{Michigan Math.\ J., {\bf 28} (1981), 97--116}







\bibitem{PTY}\listing{M.\ Passare, A.\ Tsikh, A.\ Yger}
{Residue currents of the Bochner-Martinelli type}
{Publ.\ Mat.,  {\bf 44} (2000), 85-117}







\bibitem{Sam}\listing{H.\ Samuelsson}
{Regularizations of products of residue and principal value currents}  
{J. Funct. Anal.,  {\bf 239}  (2006),  566--593}





\bibitem{Scheja}
\listing{G.\ Scheja}
{Riemannsche Hebbarkeitss\"atze f\"ur Cohomologieklassen}
{Math.\ Ann.,  {\bf 144}  (1961),  345--360}
{}{}






\bibitem{Sko}\listing{H.\ Skoda}
{Application des techniques $L\sp{2}$ \`a la th\'eorie des id\'eaux d'une alg\`ebre de 
fonctions holomorphes avec poids}{Ann.\ Sci.\   \'Ecole Norm.\ Sup., {\bf 5} (1972),
545--579}




\bibitem{Sz}\listing{J.\ Sznajdman}
{An elementary proof of the Brian\c{c}on-Skoda theorem}
{Preprint Gothenburg 2008,\ Available at  arXiv:0807.0142}{}{}{}{}




\bibitem{W}\listing{E.\ Wulcan}
{Residue currents of monomial ideals}
{Indiana Univ.\ Math.\ J.,  {\bf 56}  (2007),  365--388}






\end{thebibliography}
\end{document}